# π and prime numbers

by
Simon Plouffe
April 22, 2023

Abstract

The starting point is a function due to Dirichlet, the function $\beta(s)$

$$\beta(s) = \sum_{n=0}^{\infty} \frac{(-1)^n}{(2n+1)^s}$$

If $s = 7$ the value is $\frac{61\pi^7}{184320}$ by isolating 61 we find first of all an approximation of 61, indeed

$$61 \approx \frac{184320}{\pi^7} = 61.0271871\ldots$$

Or more precisely

$$61 = \frac{2^8\, 6!}{\pi^7}\left[1 - \frac{1}{3^7} + \frac{1}{5^7} - \frac{1}{7^7}\ldots\right]$$

So, by inverting this series we get a good approximation of a prime number and if we push further an equality. Here the number 61 is the 6th Euler number. The largest known Euler number that is prime is $E_{510}$ a number of 1062 digits, there does not seem to be any other before $E_{100000}$. We take here the numbers generated by the series $\frac{1}{\cos(x)}$ where all the $E_n$ are positive. It should also be noted that the approximation provided by this process gives us

$$E_{510} \approx \frac{2^{512}\, 510!}{\pi^{511}}$$

which is excellent since 243 exact decimals out of 1062. The goal of this article is to find all possible expressions that can represent a prime number.

# The Euler, Bernoulli and Tangent numbers

As indicated in the summary, the Taylor expansion of the function $1/\cos(x)$ function provides the Euler numbers and if some of these are prime then we will have a representation of this prime. It remains to find which other primes can be represented. The Euler-Tangent numbers for example are obtained with the expression

$$\frac{1 + \sin(x)}{\cos(x)}$$

Which once developed in exponential series generates the numbers

1, 1, 1, 2, 5, 16, 61, 272, 1385, 7936, 50521, 353792, ... (sequence A000111 of the OEIS catalog) whose general term, $a(n) \approx \frac{2^{n+2} n!}{\pi^{n+1}}$. We notice that 2 is a factor in the initial expression, we can then represent $17 = \frac{272}{16} \approx \frac{32*7!}{\pi^8}$ as well as 277 since $a(8) = 1385 \approx \frac{1024*8!}{\pi^9}$ and thus $277 \approx \frac{8257536}{\pi^9}$. We can then obtain an approximation of this prime as well as the exact expression by taking the series. Our list of primes is getting longer.

The list so far is 2, 5, 17, 31, 61, 277, 691 and $E_{510}$.

The expressions are always of the form

$$p \approx \frac{c^m k!}{d \pi^n} \tag{1}$$

Where c and d are constants, in the examples above, $c = 2$ and d can be integer or algebraic. These first few are just the tip of the iceberg for an obvious reason. It is easy to generate numbers like Euler or Bernoulli with other trigonometric expressions like

$$\frac{\sin(x) + \cos(3x)}{\cos(4x)}$$

The development in Taylor series gives the following sequence: 1, 1, 7, 47, 497, 6241, 95767, 1704527, 34741217, ...
We see 7, 47. So if we can calculate the general term a(n) of this sequence we can add the first ones found to the list. Indeed, the general term is

$$a(n) \approx \frac{8^n \Gamma(n)}{2\sqrt{4 + 2\sqrt{2}} \, \pi^n}$$

That we find in this way.

1) We develop in series the trigonometric expression

2) We collect the coefficients with n!
3) We calculate $a(n+1)/a(n)$.
4) We calculate the differences term by term to find the constant $8/\pi$.
5) We go backwards (retro-engineering) to isolate the constant d in the expression (1).

This makes it possible to obtain that
$$7 \approx \frac{8^3}{\sqrt{4+2\sqrt{2}}\,\pi^3}, \quad 47 \approx \frac{8^3\,4!}{\sqrt{4+2\sqrt{2}}\,\pi^4}$$

The approximations found are coarse when n is small but the 200<sup>th</sup> term is of the order of $10^{454}$ and $a(200)$ has a precision of 95 decimal digits.

If we look at other functions like the ones above, we find other prime numbers, here are some examples. All expressions are developed in exponential generating series. In **bold**, a term over 2 is considered. *Annnnnn* refers to the Online Encyclopedia of Integer Sequences: OEIS.org.

| Function | Sequence | General term | Prime | Example |
|---|---|---|---|---|
| $\dfrac{\sin(x)}{\cos(2x)}$ | A000464<br>1, 11, 361, 24611,<br>2873041, 512343611, | $\dfrac{(2n+1)!\,2^{4n+7/2}}{\pi^{2n+2}}$ | 11, 13, 19 | $11 \approx \dfrac{768\sqrt{2}}{\pi^4}$ |
| $\dfrac{1-2\cos(x)}{1-\sin(x)}$ | A000708<br>-1, -1, 0, 1, 6, 29, 150, 841,<br>5166, 34649, | $\dfrac{2^n\,n!}{\pi^n}$ | 3, 29, | |
| $\dfrac{\cos(x)}{\cos(2x)}$ | A000281<br>1, 3, 57, 2763, 250737, | $\dfrac{\sqrt{2}\,2^n\,(2n)!\,4^{2k}}{\pi^{2n+1}}$ | 3, 19, | $3 \approx \dfrac{64\sqrt{2}}{\pi^3}$ |
| $\dfrac{1}{\cos(x)}$ | A000364<br>1, 5, 61, 1385, | $\dfrac{2^{n+2}\,n!}{\pi^{n+1}}$ | 5, 61, 277, | $61 \approx \dfrac{2^8\,6!}{\pi^7}$ |
| $\dfrac{\sin(x)+\cos(3x)}{\cos(4x)}$ | A006873<br>1, 1, 7, 47, 497, 6241, | $\dfrac{8^n\,\Gamma(n)}{2\sqrt{4+2\sqrt{2}}\,\pi^n}$ | 7, 47, | $7 \approx \dfrac{8^3}{\sqrt{4+2\sqrt{2}}\,\pi^3}$ |
| $\dfrac{1+\sin(x)}{\cos(x)}$ | A000111<br>1, 1, 1, 2, 5, 16, 61, 272, ... | $\dfrac{2^{n+2}\,n!}{\pi^{n+1}}$ | 2, 5, 17, | $17 \approx \dfrac{32*7!}{\pi^8}$ |
| $\dfrac{\cos(x)+\cos(2x)+\cos(3x)}{\cos(x)}$ | 3, 11, 29, 191, 871, 52571, | $\dfrac{2^{2n}\,(2n-2)!}{\pi^{(2n-1)}}$ | 3, 11, 29,<br>191, 52571 | $52571 \approx \dfrac{2^{12}\,10!}{\pi^{11}}$ |

| $\dfrac{\cos(2x)+\sin(2x)}{\cos(2x)}$ | A012393 2,16,512,34816,4063232, | $\dfrac{2^{4n}\,(2n-1)!}{\pi^{2n}}$ | 17,31,691 | $691 \approx \dfrac{2^4\,11!}{\pi^{12}}$ |

We could continue these examples to find other prime numbers, these few functions are among the simplest. In addition to isolated primes if a term is composed of 2 factors, it is always possible to obtain a new identity as in the case of the 8th Euler number : 1385 which is $5 \cdot 277$. Since we already have 5, this provides an approximation for 277. The process works but it is a bit tedious to go and get any prime and illusory to find an expression for each prime. For the moment, this short list is only a sequence of exotic expressions giving almost a prime number.

## Other forms giving prime numbers

Ramanujan found this formula with Bernoulli numbers and an unusual sum.

$$\sum_{n=1}^{\infty} \frac{n^{4k+1}}{e^{2\pi n}-1} = \frac{B_{4n+2}}{2(4n+2)} \qquad (2)$$

This gives for $4k+1 = 13$ the identity

$$\sum_{n=1}^{\infty} \frac{n^{13}}{e^{2\pi n}-1} = \frac{1}{24}$$

If k is much higher, we have

$$24 \sum_{n=1}^{\infty} \frac{n^{673}}{e^{2\pi n}-1} = 1563446 \ldots 036059151$$

A prime number of 1077 digits, with the exponent 22607 we have one of 71299 digits, it is the largest known of its kind. What is remarkable for this sum is the affinity with once again the number $\pi$. Remarkable also because of the exceptional precision.

Indeed, it is not difficult to establish that in general we will have this approximation

$$\sum_{n=1}^{\infty} \frac{n^k}{e^{r\pi n}-1} \approx \frac{k!}{(r\pi)^{k+1}} \qquad (3)$$

So, if the result of the infinite sum is prime we will have another representation of these approximations. If k is high, the approximation is spectacular, here with the exponent 673 has 202 exact decimals.

These sums (2) are known, what is less known is to extend the Ramanujan formula to two terms.

Using an example with 691 we have (691 is the numerator of the 12th Bernoulli number).

$$691 = 16 \sum_{n=1}^{\infty} \frac{n^{11}}{e^{\pi n} - 1} - 2^{16} \sum_{n=1}^{\infty} \frac{n^{11}}{e^{4\pi n} - 1} \qquad (4)$$

The first sum already provides a good approximation of $691 \approx \frac{2^4 11!}{\pi^{12}}$ Again, this approximation can be made exact by adding either the missing term or a part of the series of the 2nd member of the equation. This equation has a strong resemblance to one given by Jean-Pierre Serre (A course in arithmetic p.157).

$$Eis_6 = 1 + \frac{65520}{691} \sum_{n=1}^{\infty} \sigma_{11}(n) q^n \qquad (5)$$

In disguise the series (2) are Eisenstein series where the variable $q^n$ becomes $e^{q\pi n}$. The identity (4) with 691 is just a rewriting of (5).

A summary analysis of the identity (2) and (3) shows that we have either 1 term or 2 terms which give a rational, for example

$$1 = 24 \sum_{n=1}^{\infty} \frac{n^3}{e^{\pi n} - 1} - 264 \sum_{n=1}^{\infty} \frac{n^3}{e^{2\pi n} - 1} \qquad (6)$$

We can obtain an infinity of them as (6) if the exponent is of the form $4k - 1$. The coefficients grow with k, when k = 5 we have

$$221930581 = 16 \sum_{n=1}^{\infty} \frac{n^{19}}{e^{\pi n} - 1} - 2^{24} \sum_{n=1}^{\infty} \frac{n^{19}}{e^{4\pi n} - 1}$$

$$221930581 \approx \frac{1946319750384384375}{\pi^{20}} = 221930369.2868 \ldots \qquad (7)$$

But $221930581 = 31 \cdot 41 \cdot 283 \cdot 617$, is not prime. Again, these examples also fall into the category of isolated exotic expressions. The approximations are very precise but there are too few examples. So, we are limited to the few cases of primes appearing in these 2-term sums. But there is a trick, by taking the sum or the difference of 2 different values (different k). For example

$$31 = -504 \sum_{n=1}^{\infty} \frac{n^5}{e^{\pi n} - 1} - 32256 \sum_{n=1}^{\infty} \frac{n^5}{e^{4\pi n} - 1} + 64 \sum_{n=1}^{\infty} \frac{n^7}{e^{\pi n} - 1} - 2^{14} \sum_{n=1}^{\infty} \frac{n^7}{e^{4\pi n} - 1}$$

Who provides the approximation

$$31 \approx \frac{61425}{\pi^6} - \frac{321300}{\pi^8} = 30.0299856 \ldots$$

The approximations are less spectacular here, but what we gain is generality. By combining 4 terms of the type and with r = 1 or 4 and k odd.

$$S(k,r) = \sum_{n=1}^{\infty} \frac{n^k}{e^{r\pi n} - 1}$$

one can easily obtain several representations for each prime even from 3. 2 having been obtained with Tangent numbers. Not only does it work quite well for all the primes up to 1000000 but also one can choose an arbitrary prime.

A test has been done with the first ones of the form $10^m + c_m$ with m = 1 up to 101. These are the first prime numbers greater than $10^m$ if m = 100 then $c_m = 267$. It is the sequence 2, 11, 101, 1009, ... (the sequence A003617 of the OEIS catalog).

Example with m = 18,

$$10^{18} + 3 \approx -\frac{57780097674288877833353934375}{8\pi^{20}} + \frac{124624995799434279343497115 3125}{8\pi^{24}}$$

Which is 1000000676938336801.2703…

The true value of the former is

$$10^{18} + 3 = 5937369506 \sum_{n=1}^{\infty} \frac{n^{19}}{e^{\pi n} - 1} - 6225783167123456 \sum_{n=1}^{\infty} \frac{n^{19}}{e^{4\pi n} - 1}$$
$$- 6025884 \sum_{n=1}^{\infty} \frac{n^{23}}{e^{\pi n} - 1} + 101097557458944 \sum_{n=1}^{\infty} \frac{n^{23}}{e^{4\pi n} - 1}$$

And for each of the terms of this sum we have the approximation (3), obviously the coefficients of the sums are prime to each other, two by two.

$$\frac{k!}{(r\pi)^{k+1}} \approx S(k,r)$$

which can be corrected with $\zeta(k + 1)$ to give an even more precise rational value if k is large, in other words

$$S(k,r) \approx \frac{k!}{(r\pi)^{k+1}} \zeta(k+1) \quad (8)$$

which is rational when k is odd. More exactly when $k = 4n + 1$. When $k = 4n + 3$ it takes 2 terms in $S(k, r)$. If we check the previous calculation we obtain the sum of 4 rationals.

$$10^{18} + 3 = -\frac{67943139037162110976}{825} + \frac{518365013406083}{6600}$$
$$+ \frac{248914919490662170624}{1365} - \frac{118691882844287}{10920}$$

There is something rather mysterious going on between equation (8) which gives a rational and the true value of $S(k,r)$ which is irrational except if $r = 2$ and k is of the form 4n+3 (Ramanujan formula). Here what we found is that the formula (8) gives the true value if we take 2 terms in $S(k,r)$. Here is an example,

We know that the following identity is true:

$$1 = 16 \sum_{n=1}^{\infty} \frac{n^3}{e^{\pi n} - 1} - 256 \sum_{n=1}^{\infty} \frac{n^3}{e^{4\pi n} - 1}$$

But knowing that $S(3,1)$ and $S(3,4)$ are likely to be irrational taken separately. Our approximate formula (8) gives us $S(3,1) \approx \frac{1}{15}$ and $S(3,4) \approx \frac{1}{3840}$ and if we add the coefficients (here 16 and -256), we get the true value of 1. Indeed :

$$16\left(\frac{1}{15}\right) - 256\left(\frac{1}{3840}\right) = 1$$

We can see that when k is small, here 3, the approximation is not very good and despite that, the equation becomes exact if we sum 2 terms. If k is larger, the approximation is much better and the result of the sum of 2 terms is again exact as with the case of $p = 10^{18} + 3$ and $k = 19$ and 23.

So, with a little program using to find the coefficients (lindep of pari-GP or PSLQ in Maple) of the sums $S(k,r)$ for each prime, we can easily reach 1000000 in a few minutes. As the size of the prime increases, there are more exact candidate sums and as many approximations.

Given that the choice of the prime to represent seems arbitrary it is reasonable to conjecture that for any prime P we have the sum

$$P = a \sum_{n=1}^{\infty} \frac{n^{k_1}}{e^{\pi n} - 1} + b \sum_{n=1}^{\infty} \frac{n^{k_1}}{e^{4\pi n} - 1} + c \sum_{n=1}^{\infty} \frac{n^{k_2}}{e^{\pi n} - 1} + d \sum_{n=1}^{\infty} \frac{n^{k_2}}{e^{4\pi n} - 1} \quad (7)$$

Here $k_1$ and $k_2$ are odd and different, a, b, c, and d are integers. Moreover, the prime P is approximated by

$$P \approx u \frac{k_1!}{(r\pi)^{k_1+1}} + v \frac{k_2!}{(r\pi)^{k_2+1}}$$

Where u and v are integers too. It remains to know what $S(k,r)$ represents exactly. From the Ramanujan identity we can easily find those with r = 1 or 4. Bill Gosper had already explored some of these identities, including this one

$$\sum_{n=1}^{\infty} \frac{n^3}{e^{2\pi n/7} - 1} = \frac{-1}{240} + \frac{1}{320}\left(301 + 210\sqrt{2}\,7^{1/4} + 120\sqrt{7} + 90\sqrt{2}\,7^{3/4}\right)\frac{\pi^2}{\Gamma(\frac{3}{4})^8}$$

This allows us to find the explicit values for r = 1, 2, 4 in $S(k,r)$.

| $S(k,r) = \sum_{n=1}^{\infty} \frac{n^k}{e^{r\pi n}-1}$ | Value | Approximation |
|---|---|---|
| $S(3,1)$ | $\frac{11}{320}\frac{\pi^2}{\Gamma(3/4)^8} - \frac{1}{240}$ | $\frac{6}{\pi^4}$ |
| $S(3,2)$ | $\frac{1}{320}\frac{\pi^2}{\Gamma(3/4)^8} - \frac{1}{240}$ | $\frac{3}{8\,\pi^4}$ |
| $S(3,4)$ | $\frac{11}{5120}\frac{\pi^2}{\Gamma(3/4)^8} - \frac{1}{240}$ | $\frac{3}{128\,\pi^4}$ |
| $S(5,1)$ | $\frac{3}{64}\frac{\pi^3}{\Gamma(3/4)^{12}} + \frac{1}{504}$ | $\frac{120}{\pi^6}$ |
| $S(5,2)$ | $\frac{1}{504}$ | $\frac{15}{8\,\pi^6}$ |
| $S(5,4)$ | $\frac{-3}{2^{12}}\frac{\pi^3}{\Gamma(3/4)^{12}} + \frac{1}{504}$ | $\frac{15}{2^9\,\pi^6}$ |
| $S(7,1)$ | $\frac{363}{2^{12}\,5}\frac{\pi^4}{\Gamma(3/4)^{16}} - \frac{1}{480}$ | $\frac{7!}{\pi^8}$ |
| $S(7,2)$ | $\frac{3}{2^{12}\,5}\frac{\pi^4}{\Gamma(3/4)^{16}} - \frac{1}{480}$ | $\frac{315}{2^4\,\pi^8}$ |
| $S(7,4)$ | $\frac{363}{2^{17}\,5}\frac{\pi^4}{\Gamma(3/4)^{16}} - \frac{1}{480}$ | $\frac{315}{2^{12}\,\pi^6}$ |
| $S(9,1)$ | $\frac{189}{2^8}\frac{\pi^5}{\Gamma(3/4)^{20}} + \frac{1}{264}$ | $\frac{9!}{\pi^{10}}$ |
| $S(9,2)$ | $\frac{1}{264}$ | $\frac{2835}{8\,\pi^{10}}$ |
| $S(9,4)$ | $\frac{189}{2^{18}}\frac{\pi^5}{\Gamma(3/4)^{20}} + \frac{1}{264}$ | $\frac{2835}{2^{13}\,\pi^{10}}$ |
| $S(11,1)$ | $\frac{393309}{66560}\frac{\pi^6}{\Gamma(3/4)^{24}} - \frac{691}{65520}$ | $\frac{11!}{\pi^{12}}$ |
| $S(11,2)$ | $\frac{189}{66560}\frac{\pi^6}{\Gamma(3/4)^{24}} - \frac{691}{65520}$ | $\frac{155925}{2^4\,\pi^{12}}$ |
| $S(11,4)$ | $\frac{393309}{2^{22}\,5\,13}\frac{\pi^6}{\Gamma(3/4)^{24}} - \frac{691}{65520}$ | $\frac{155925}{2^{16}\,\pi^{12}}$ |
| $S(13,1)$ | $\frac{68607}{2^{10}}\frac{\pi^7}{\Gamma(3/4)^{28}} - \frac{1}{24}$ | $\frac{13!}{\pi^{14}}$ |
| $S(13,2)$ | $\dfrac{1}{24}$ | $\dfrac{6081075}{2^4\,\pi^{14}}$ |
| $S(13,4)$ | $\dfrac{68607}{2^{24}}\dfrac{\pi^7}{\Gamma(3/4)^{28}}+\dfrac{1}{24}$ | $\dfrac{6081075}{2^{18}\,\pi^{14}}$ |

# Appendix

Table of approximations of primes from 2.

| Prime | Approximation | Value |
|---|---|---|
| 2 | $\dfrac{196}{\pi^4}$ | 2.0121 |
| 103 | $\dfrac{184275}{8\,\pi^6} + \dfrac{2972025}{4\,\pi^8}$ | 102.2651 |
| 163 | $\dfrac{798525}{8\,\pi^6} + \dfrac{2168775}{4\,\pi^8}$ | 160.9663 |
| $10^{18} + 3$ | $\approx -\dfrac{57780097674288877833539343 75}{8\,\pi^{20}} + \dfrac{124624995799434279343497115 3125}{8\,\pi^{24}}$ | 1000000676938336801.2703 |

# Raw table

$Prime, S(k,1), S(k,4)$ [$vector\ that\ makes\ the\ sum = 0$] $Approximation$

Example

$$7 = -228\,S(5,1) + 58368\,S(5,4) + 33\,S(7,1) + 33792\,S(7,4)$$

$$7 \approx \dfrac{4578525}{84\,\pi^8} + \dfrac{95893875}{8\,\pi^{10}}$$

```
7,   7,  9,  [-1, -228,  58368,  33,  33792],  -4578525/4/Pi ^8+95893875/8/Pi ^10
23,  5,  7,  [-1,   63,   4032,  28,  -7168],  61425/8/Pi ^6+562275/4/Pi ^8
29,  5,  7,  [-1,  126,   8064,  24,  -6144] ,  61425/4/Pi ^6+240975/2/Pi ^8
31,  5,  7,  [-1,  504,  32256, -64,  16384],  61425/Pi ^6-321300/Pi ^8
31,  7,  9,  [-1, -424, 108544,  66,  67584]  -4257225/2/Pi ^8+95893875/4/Pi ^10
41,  5,  7,  [-1,  252,  16128,  16,  -4096],  61425/2/Pi ^6+80325/Pi ^8
47,  5,  7,  [-1,  315,  20160,  12,  -3072],  307125/8/Pi ^6+240975/4/Pi ^8
53,  5,  7,  [-1,  378,  24192,   8,  -2048],  184275/4/Pi ^6+80325/2/Pi ^8
59,  5,  7,  [-1,  441,  28224,   4,  -1024],  429975/8/Pi ^6+80325/4/Pi ^8
71,  5, 11,  [-1, 1890, 120960,  -4,  16384],  921375/4/Pi ^6-638512875/4/Pi ^12
71,  5,  7,  [-1,  567,  36288,  -4,   1024],  552825/8/Pi ^6-80325/4/Pi ^8
89,  5,  7,  [-1,  756,  48384, -16,   4096],  184275/2/Pi ^6-80325/Pi ^8
89,  7,  9,  [-1, -556, 142336,  99, 101376],  -11165175/4/Pi ^8+287681625/8/Pi ^10
97,  3,  9,  [-1, -498,   7968,  33,  33792],  -11205/4/Pi ^4+95893875/8/Pi ^10
97,  5,  7,  [-1,  126,   8064, 152, -38912],  61425/4/Pi ^6+1526175/2/Pi ^8
103, 5,  7,  [-1,  189,  12096, 148, -37888],  184275/8/Pi ^6+2972025/4/Pi ^8
109, 5,  7,  [-1,  252,  16128, 144, -36864],  61425/2/Pi ^6+722925/Pi ^8
109, 7,  9,  [-1,  -36,   9216,  33,  33792],  -722925/4/Pi ^8+95893875/8/Pi ^10
113, 3,  9,  [-1, -242,   3872,  33,  33792],  -5445/4/Pi ^4+95893875/8/Pi ^10
113, 5,  7,  [-1, 1008,  64512, -32,   8192],  122850/Pi ^6-160650/Pi ^8
127, 5, 11,  [-1,  315,  20160,   2,  -8192],  307125/8/Pi ^6+638512875/8/Pi ^12
127, 5,  7,  [-1,  441,  28224, 132, -33792],  429975/8/Pi ^6+2650725/4/Pi ^8
131, 3,  9,  [-1,   46,   -736,  33,  33792],  1035/4/Pi ^4+95893875/8/Pi ^10
131, 5,  7,  [-1,  126,   8064, 216, -55296],  61425/4/Pi ^6+2168775/2/Pi ^8
137, 3,  9,  [-1,  142,  -2272,  33,  33792],  3195/4/Pi ^4+95893875/8/Pi ^10
137, 5,  7,  [-1,  189,  12096, 212, -54272],  184275/8/Pi ^6+4257225/4/Pi ^8
139, 5, 11,  [-1, 3087, 197568,  -6,  24576],  3009825/8/Pi ^6-1915538625/8/Pi ^12
139, 5,  7,  [-1,  567,  36288, 124, -31744],  552825/8/Pi ^6+2490075/4/Pi ^8
149, 3,  9,  [-1,  334,  -5344,  33,  33792],  7515/4/Pi ^4+95893875/8/Pi ^10
149, 5,  7,  [-1,  315,  20160, 204, -52224],  307125/8/Pi ^6+4096575/4/Pi ^8
151, 3,  9,  [-1,  366,  -5856,  33,  33792],  8235/4/Pi ^4+95893875/8/Pi ^10
151, 5,  7,  [-1,  693,  44352, 116, -29696],  675675/8/Pi ^6+2329425/4/Pi ^8
```

```
157, 3, 9, [-1, 462, -7392, 33, 33792], 10395/4/Pi ^4+95893875/8/Pi ^10
157, 5, 7, [-1, 756, 48384, 112, -28672], 184275/2/Pi ^6+562275/Pi ^8
157, 7, 9, [-1, -428, 109568, 99, 101376], -8594775/4/Pi ^8+287681625/8/Pi ^10
163, 3, 9, [-1, 558, -8928, 33, 33792], 12555/4/Pi ^4+95893875/8/Pi ^10
163, 5, 7, [-1, 819, 52416, 108, -27648], 798525/8/Pi ^6+2168775/4/Pi ^8
167, 7, 9, [-1, -168, 43008, 66, 67584], -1686825/2/Pi ^8+95893875/4/Pi ^10
167, 9, 11, [-1, 132, 135168, -8, 32768], 95893875/2/Pi ^10-638512875/2/Pi ^12
173, 3, 9, [-1, 718, -11488, 33, 33792], 16155/4/Pi ^4+95893875/8/Pi ^10
179, 3, 9, [-1, 814, -13024, 33, 33792], 18315/4/Pi ^4+95893875/8/Pi ^10
181, 3, 9, [-1, 846, -13536, 33, 33792], 19035/4/Pi ^4+95893875/8/Pi ^10
181, 5, 7, [-1, 1008, 64512, 96, -24576], 122850/Pi ^6+481950/Pi ^8
191, 5, 7, [-1, 756, 48384, 176, -45056], 184275/2/Pi ^6+883575/Pi ^8
```